\documentclass[letterpaper, 10 pt, conference]{ieeeconf}  
\usepackage[top=0.76in, bottom = 0.8in, left = 0.77in, right=0.77in]{geometry}
\usepackage{amsmath,amssymb,amsfonts}
\usepackage{tikz}
\usetikzlibrary{shapes,arrows}
\usepackage{verbatim}

\DeclareMathOperator{\diag}{diag}
\usepackage{graphicx}
\usepackage{etoolbox}
\usepackage{cases}
\usepackage{algorithmic}
\usepackage{graphicx}
\usepackage{textcomp}
\usepackage{xcolor}
\usepackage{array}

\usepackage{cite}
\usepackage{mathrsfs}
\newtheorem{theorem}{Theorem}
\newtheorem{lemma}{Lemma}[section]
\newtheorem{remark}{Remark}
\newtheorem{definition}{Definition}

\newtheorem{corollary}{Corollary}[theorem]

\newtheorem{condition}{Condition}
\usepackage{comment}
\usepackage{tabularx}
\usepackage{multirow}
\IEEEoverridecommandlockouts                              

\title{\LARGE \bf
A Fundamental Convergence Rate Bound for Gradient Based Online Optimization Algorithms with Exact Tracking}

\author{Alex (Xinting) Wu, Ian R. Petersen, Iman Shames
\thanks{*This work was supported by the Australian Research Council under grants DP230102443 and DP210102454.
}
\thanks{A. Wu, I. R. Petersen and I. Shames are with the CIICADA Lab, School of Engineering, The Australian National University, Canberra, ACT 2601,
Australia (e-mail: u5847417@anu.edu.au; i.r.petersen@gmail.com; iman.shames@anu.edu.au).}%
}%

\begin{document}

\maketitle
\thispagestyle{empty}
\pagestyle{empty}

\begin{abstract}
In this paper, we consider algorithms with integral action for solving online optimization problems characterized by quadratic cost functions with a time-varying optimal point described by an $(n-1)$th order polynomial. Using a version of the internal model principle, the optimization algorithms under consideration are required to incorporate a discrete time $n$-th order integrator in order to achieve exact tracking. By using results on an optimal gain margin problem, we obtain a fundamental convergence rate bound for the class of linear gradient based algorithms exactly tracking a time-varying optimal point. This convergence rate bound is given by $ \left(\frac{\sqrt{\kappa} - 1 }{\sqrt{\kappa} + 1}\right)^{\frac{1}{n}}$, where $\kappa$ is the condition number for the set of cost functions under consideration. Using our approach, we also construct algorithms which achieve the optimal convergence rate as well as zero steady-state error when tracking a time-varying optimal point.

\end{abstract}

\section{INTRODUCTION}
Online optimization \cite{Introonline,Zampieri2023,madden2021bounds} has emerged as an important area of research with implications across various fields, including real-time control systems, signal processing, and machine learning. An online optimization problem involves making a sequence of decisions in real time, where the objective function changes continually. The primary challenge lies in efficiently solving these problems while ensuring that the time-varying solutions converge asymptotically. This motivates the development of algorithms capable of fast convergence and zero steady-state error. However, many existing methods for online optimization such as those of \cite{Introonline,madden2021bounds}, do not achieve zero steady-state error. 

In recent years, various methods have been employed for designing optimization algorithms and analyzing their convergence behavior, utilizing techniques from control theory, including the internal model principle \cite{Zampieri2023}, integral quadratic constraints (IQC) \cite{doi:10.1137/15M1009597, IQCScherer}, Lyapunov-based analysis \cite{Lessard2023L} and dissipative systems \cite{lessard2022}. More recently, \cite{UGRINOVSKII2023111129, Jovanovic2024Tann, frequencyTannen2024} build on the connection between the optimal gain margin problem \cite{Tannenbaum1985} and optimization algorithms, with a particular focus on the heavy ball method. This link has provided new insights into understanding the convergence behavior and stability of optimization algorithms by drawing parallels with control theory concepts.  

In this paper, we consider a class of online optimization problems with quadratic cost functions where the optimal point varies polynomially with the time step. The internal model principle suggests that in order to achieve zero steady-state error, an algorithm needs to incorporate a suitable number of discrete time integrators \cite{AlexACCTracking}. The paper extends the results of \cite{AlexACCTracking} by using the optimal gain margin approach of \cite{Tannenbaum1985} to establish a fundamental convergence rate bound of $ \left(\frac{\sqrt{\kappa} - 1 }{\sqrt{\kappa} + 1}\right)^{\frac{1}{n}}$ for all linear gradient based algorithms which achieve zero steady-state error when tracking a polynomially varying optimal point of order $n-1$. Here, $\kappa$ is the condition number for the class of cost functions, which will be defined in the sequel. In addition, for any set of cost functions which includes the set of quadratic cost functions under consideration, will have a convergence rate which satisfies this bound. Moreover, the paper shows that this convergence rate bound is tight in that we give a procedure to construct algorithms which achieve the bound.

A conference version of this paper was presented in \cite{AlexECC2025}. In this paper, we extend the results of \cite{AlexECC2025} by generalizing from a linearly varying optimal point to a polynomially varying optimal point and by establishing a fundamental convergence rate bound for gradient based optimization algorithms with exact tracking. 

This paper is organized as follows: Section \ref{sec2} presents the problem formulation and preliminary results. In Section \ref{sec3}, we describe and prove our main results. Section \ref{IlluExam} presents an illustrative example. Finally, the paper is concluded in Section \ref{sec5}.

\section{Problem formulation and Preliminary Results} \label{sec2}
In this section, we formulate the class of optimization problems under consideration and review some existing optimization methods and their convergence rates, alongside relevant preliminary results.

\subsection{Problem formulation}
We consider unconstrained optimization problems of the form:
\begin{equation}
\label{opproblem}
    x^*(t) = \arg \min f(x,t),
\end{equation}
where $f: \mathbb{R}^{p+1} \rightarrow \mathbb{R}$ is a quadratic cost function with a unique minimum attained at $x^*(t)$. 

To address this problem, we consider gradient based algorithms similar to the general algorithms considered in \cite{UGRINOVSKII2023111129} of the following form:
\begin{align}
    x(t+1) &= x(t) +\sum_{j=0}^{k-1} \beta_j (x(t-j)-x(t-j-1)) \nonumber\\& \quad -\sum_{j=0}^{k}\alpha_j \nabla_x f(x(t-j),t-j), \label{uprule}
\end{align}
where $k \in \mathbb{N}, \alpha_j \in \mathbb{R}, \beta_j \in \mathbb{R}$ and
\begin{equation}
    \sum_{j=0}^{k}\alpha_j \neq 0. \label{alpha0}\\
\end{equation}

In (\ref{uprule}), $x(t) \in\mathbb{R}^p$ represents the current iteration point and $t \in {0, 1, \ldots}$ denotes the iteration index. The number $k$ denotes the number of past iterates considered in the optimization algorithm and the number of gradient evaluations at each iteration. Since $k$ corresponds to the number of past iterates, the algorithm needs be initialized with the values $x(0), x(1), \dots, x(k)$.
The parameters $\alpha_j$ and $\beta_j$ are the algorithm parameters, and $\nabla_x f(x(t),t)$ denotes the gradient of the cost function. It is straightforward to verify that this class of algorithms includes all linear gradient based optimization algorithms with constant step size.

We consider a class of time-varying quadratic cost functions defined as follows.
\begin{definition}
    Given $L \geq m >0$, let $\mathscr{F}^q_{m,L}$ denote the class of time-varying quadratic cost functions $f(x,t)$ of the following form:
    \begin{align}
    \label{quadrafunc}
        f(x,t) &= \frac{1}{2}\left(x(t)-x^*(t) \right)^T \Delta \left(x(t)-x^*(t) \right) \nonumber \\
        &\quad +c( x^*(t)),
    \end{align}
    where the time-varying optimal point $x^*(t)$ is of the form
    \begin{equation}
        \label{x*t}
        x^*(t) = a_0+a_1t+..+a_{n-1}t^{n-1},
    \end{equation}
    and $a_i\in \mathbb{R}^p$ are constant vectors. 
    In (\ref{quadrafunc}),
    \begin{equation}
    \label{classfunc}
        mI \leq \Delta \leq LI, \quad \Delta = \Delta^T \in \mathbb{R}^{p \times p}, c( x^*(t))\in\mathbb{R}.
    \end{equation} 
    We also define $\kappa  = \frac{L}{m}$.
\end{definition}

In the sequel, we will show that in order for the algorithm (\ref{uprule}) to track the polynomially varying optimal point defined in (\ref{x*t}), it is necessary that it satisfies the following condition.

\begin{condition}
    \label{assumpbeta}
    The algorithm (\ref{uprule}) is such that 
    \begin{equation}
        \sum_{j=0}^{k-1} \beta_j(\hat{k}-j-1)_{r}=(\hat{k})_{r}, \label{betak-1}
    \end{equation}
    for any $0\leq r \leq n-2$ and any $\hat{k} \geq k$,
    where $(\cdot)_{r}$ denotes the falling factorial of order $r$ \cite{concretemath}; i.e., $(\hat{k})_{r} = \hat{k}(\hat{k}-1)...(\hat{k}-r+1)$.
\end{condition}

In the sequel, we will demonstrate that Condition \ref{assumpbeta} corresponds to the requirement that the algorithm includes at least $n$ discrete time integrators and so obtains zero steady-state error with an optimal point of the form (\ref{x*t}) according to the internal model principle; e.g., see \cite{Zampieri2023}.

In this paper, we consider online optimization algorithms of the form (\ref{uprule}) to track the minimum in (\ref{uprule}) where $f(x,t) \in \mathscr{F}^q_{m,L}$. The algorithms (\ref{uprule}) can be reformulated as uncertain linear feedback systems, enabling us to analyze their convergence behavior through the lens of an optimal gain margin problem. 
This work extends the findings of \cite{UGRINOVSKII2023111129} by considering online optimization problems in which the optimal point varies polynomially with time. The paper extend the findings of \cite{AlexACCTracking} by considering the optimal convergence rate with respect to all algorithms of the form (\ref{uprule}).
In particular, we establish a fundamental convergence rate bound for all algorithms of the form (\ref{uprule}) that achieve zero steady-state error for time-varying cost functions of the form $f(x,t) \in \mathscr{F}^q_{m,L}$. This convergence rate bound is given by the formula $\left(\frac{\sqrt{\kappa} - 1 }{\sqrt{\kappa} + 1}\right)^{\frac{1}{n}}=\left(\frac{\sqrt{L} - \sqrt{m}}{\sqrt{L} + \sqrt{m}}\right)^{\frac{1}{n}}$. In addition, we show that this bound is tight in the sense that we can construct an algorithm which achieves this bound. This algorithm is constructed using a Blaschke product solution to a corresponding Nevanlinna Pick interpolation problem \cite{pick1915,Nevanlinna1919}. 

\subsection{The Heavy Ball Method}
The heavy ball method \cite{POLYAK19641} is one of the most well-known methods for solving a time-invariant version of optimization problem (\ref{opproblem}) in which the cost function $f(x,t)$ is independent of $t$. This approach enhances the standard gradient descent algorithm by incorporating a momentum term, resulting in an optimization algorithm of the form:
    \begin{equation} 
        \label{HB} 
        x(t+1) = x(t) - \alpha \nabla f(x(t)) + \beta(x(t) - x(t-1)). 
    \end{equation}

The heavy ball method accelerates convergence by combining the current gradient with information from the previous step. With Polyak's choice of parameters, this method exhibits the fastest convergence rate for time-invariant cost functions $f(x) \in \mathscr{F}^q_{m,L}$ \cite{UGRINOVSKII2023111129}. These parameters are defined as:
\begin{equation}
\label{polypara}
    \begin{aligned}
        \alpha = \alpha_{\text{\scriptsize HB}} =&\frac{4}{(\sqrt{L}+\sqrt{m})^2}= \frac{4}{m(1+\sqrt{\kappa})^2},\\
        \beta = \beta_{\text{\scriptsize HB}} =& \frac{(\sqrt{L}-\sqrt{m})^2}{(\sqrt{L}+\sqrt{m})^2}= \frac{(\sqrt{\kappa}-1)^2}{(\sqrt{\kappa}+1)^2}.
    \end{aligned}
\end{equation}
The paper \cite{UGRINOVSKII2023111129} has shown that, for any time-invariant cost functions $f(x) \in \mathscr{F}^q_{m,L}$, the heavy ball method yields the optimal worst-case asymptotic convergence rate, among all algorithms characterized by the following general form:
\begin{align}
    x(t+1) &= x(t) +\sum_{j=0}^{k-1} \beta_j (x(t-j)-x(t-j-1)) \nonumber\\& \quad -\sum_{j=0}^{l}\alpha_j \nabla_x f(y(t-j)),\nonumber\\
    y(t) &= \sum_{\nu=0}^{k-l} \gamma_\nu x(t-\nu), \quad t = 0,1,\dots \label{methodV}
\end{align}
where 
\begin{equation*}
    \sum_{j=0}^l \alpha_j \neq 0, \quad \sum_{\nu=0}^{k-l}\gamma_\nu=1.
\end{equation*}
The convergence rate of the heavy ball method with parameters (\ref{polypara}) is calculated as in \cite{UGRINOVSKII2023111129,POLYAK19641}:
\begin{equation}
    \label{rateHB}
    \rho_{\text{\scriptsize HB}}=\frac{\sqrt{\kappa} - 1 }{\sqrt{\kappa} + 1} = \frac{\sqrt{L}-\sqrt{m}}{\sqrt{L}+\sqrt{m}}.
\end{equation}
This formula then defines a bound on the convergence rate for all algorithms of the form (\ref{methodV}) when applied to time-invariant cost functions; see \cite{UGRINOVSKII2023111129,LecOnCOptNes2018}.


\subsection{Preliminary results}
Let $\Sigma(\alpha,\beta,m,L)$ denote an algorithm corresponding to the recursion (\ref{uprule}) depending on parameters $\alpha = (\alpha_0,\dots,\alpha_k), \beta=(\beta_0,\dots,\beta_{k-1})$ and the set of cost functions $\mathscr{F}^q_{m,L}$ defined by (\ref{quadrafunc}), (\ref{x*t}), (\ref{classfunc}) where $L>m>0$. 

The following lemma provides a time-invariant representation of the iterative algorithm (\ref{uprule}) in the case that Condition \ref{assumpbeta} is satisfied. This provides an alternative approach to the internal model principle for exact tracking of polynomially varying signals.
\begin{lemma}
    Given $\alpha$, $\beta$, $m$ and $L$, such that $L>m>0$, and an algorithm $\Sigma(\alpha,\beta,m,L)$ of the form (\ref{uprule}) such that Condition \ref{assumpbeta} is satisfied. Then (\ref{uprule}) can be rewritten as: 
    \begin{align}
    \label{upruletilx}
        \tilde{x}(t+1)&=\tilde{x}(t) +\sum_{j=0}^{k-1} \beta_j (\tilde{x}(t-j)-\tilde{x}(t-j-1))\nonumber \\
        &\quad-\sum_{j=0}^{k}\alpha_j \nabla_x \tilde{f}(\tilde{x}(t-j)),
    \end{align}
    where
    \begin{align} 
        \label{xttil}
        \tilde{x}(t) &= x(t)-x^*(t) \nonumber \\
        &= x(t) - a_0 - a_1t- a_2t^2-\dots - a_{n-1}t^{n-1}.
    \end{align}
\end{lemma}

\textit{Proof:}
To prove this lemma, note that it follows from (\ref{xttil}) that
\begin{align}
    \label{ftoftil}
    f(x,t) &= \tilde{f} \left(x(t)-x^*(t) \right)\nonumber \\ &=  \tilde{f}(\tilde{x}(t))
\end{align}

where $\tilde{f}(\tilde{x}) = \frac{1}{2} \tilde{x}^T \Delta \tilde{x}$. Also,
\begin{align}
\label{tilnabla}
    &\nabla_x f(x,t)\nonumber \\ & = \nabla_x \tilde{f}\left(x(t)-x^*(t) \right) \nonumber \\ &= \nabla_{\tilde{x}} \tilde{f}(\tilde{x}(t)).
\end{align}
It follows from (\ref{tilnabla}) that
\begin{equation}
\label{tilnabalpha}
    \sum_{j=0}^{k}\alpha_j \nabla_x f(x(t-j),t-j) = \sum_{j=0}^{k}\alpha_j \nabla_x \tilde{f}(\tilde{x}(t-j)).
\end{equation}
In addition,
\begin{align}
\label{tilsumbeta}
    &\sum_{j=0}^{k-1} \beta_j (\tilde{x}(t-j)-\tilde{x}(t-j-1)) \nonumber
    \\ & \quad= \sum_{j=0}^{k-1} \beta_j (x(t-j)-x^*(t-j)) \nonumber \\&\quad \quad-(x(t-j-1)-x^*(t-j-1)) \nonumber \\
    &\quad =  \sum_{j=0}^{k-1} \beta_j(x(t-j)-x(t-j-1)) \nonumber \\&\quad \quad-\sum_{j=0}^{k-1} \beta_j(x^*(t-j)-x^*(t-j-1)).
\end{align}

Also, it follows from (\ref{xttil}) that
\begin{equation}
    \label{xt+1til}
    \tilde{x}(t+1) = x(t+1)-x^*(t+1).
\end{equation}
Substituting (\ref{uprule}) into (\ref{xt+1til}), we obtain
\begin{align}
\label{xtil+1}
    \tilde{x}(t+1) &= x(t)-x^*(t+1) \nonumber \\&\quad+\sum_{j=0}^{k-1} \beta_j (x(t-j)-x(t-j-1)) \nonumber\\& \quad -\sum_{j=0}^{k}\alpha_j \nabla_x f(x(t-j),t-j).
\end{align}

Substituting (\ref{xttil}), (\ref{tilnabalpha}) and (\ref{tilsumbeta}) into (\ref{xtil+1}) yields
\begin{align}
\label{upruletilxf0}
    \tilde{x}(t+1)&=\tilde{x}(t) +\sum_{j=0}^{k-1} \beta_j (\tilde{x}(t-j)-\tilde{x}(t-j-1))\nonumber \\
    &\quad-\sum_{j=0}^{k}\alpha_j \nabla_x \tilde{f}(\tilde{x}(t-j)) + f_0(t),
\end{align}
where
\begin{equation}
    f_0(t)=\sum_{j=0}^{k-1} \beta_j(x^*(t-j)-x^*(t-j-1))-(x^*(t+1)-x^*(t)).
\end{equation}
Thus, we need to show that $f_0(t) \equiv 0$. Indeed, it follows from (\ref{x*t}) that $f_0(t) \equiv 0$ if 
\begin{equation}
    \label{t-jt+1}
    \sum_{j=0}^{k-1} \beta_j \left((t-j)^{\hat{n}}-(t-j-1)^{\hat{n}}\right)=(t+1)^{\hat{n}}-t^{\hat{n}},
\end{equation}
for $\hat{n}=1, 2, \dots, n-1$.
To verify (\ref{t-jt+1}), we introduce Stirling numbers of the second kind; e.g., see \cite[Section 6.1]{concretemath} which gives the following relations
\begin{equation}
    t^{\hat{n}} =\sum_{s=0}^{\hat{n}} \genfrac\{\}{0pt}{}{\hat{n}}{s} (t)_{s}
\end{equation}
where $\{ \cdot \}$ denotes the Stirling numbers of the second kind and $(t)_{s}$ denotes the falling factorial.
Therefore
\begin{align}
    \label{t+1-t}
    (t+1)^{\hat{n}}-t^{\hat{n}} &= \sum_{r=0}^{\hat{n}} \genfrac\{\}{0pt}{}{\hat{n}}{s} (t+1)_{s}-\sum_{s=0}^{\hat{n}} \genfrac\{\}{0pt}{}{\hat{n}}{s} (t)_{s} \nonumber \\
    &=\sum_{s=0}^{\hat{n}} \genfrac\{\}{0pt}{}{\hat{n}}{s}s t(t-1)(t-2)...(t-s+2) \nonumber \\
    &=\sum_{s=1}^{\hat{n}} \genfrac\{\}{0pt}{}{\hat{n}}{s}s (t)_{s-1}.
\end{align}

Similarly,
\begin{align}
    \label{strt-j}
    &(t-j)^{\hat{n}}-(t-j-1)^{\hat{n}} \nonumber\\ &\quad=\sum_{s=1}^{\hat{n}} \genfrac\{\}{0pt}{}{\hat{n}}{s}s (t-j-1)...(t-j-s+1) \nonumber \\&\quad=\sum_{s=1}^{\hat{n}} \genfrac\{\}{0pt}{}{\hat{n}}{s}s (t-j-1)_{s-1}.
\end{align}

It is straightforward to verify that $r = s-1$ and it follows from (\ref{betak-1}) that
\begin{align}
    \label{t-j-1=t}
    &\sum_{j=0}^{t-1} \beta_j(t-j-1)_{s-1} \nonumber\\ & \quad =\sum_{j=0}^{k-1} \beta_j(t-j-1)_{s-1}+\sum_{j=k}^{t-1} \beta_j(t-j-1)_{s-1}
\end{align}
Since $\beta_j = 0$ for $j>k$, it follows that
\begin{align}
    \sum_{j=0}^{t-1} \beta_j(t-j-1)_{s-1}& =\sum_{j=0}^{k-1} \beta_j(t-j-1)_{s-1} \nonumber\\ &=(t)_{s-1}.
\end{align}
Substituting (\ref{t-j-1=t}) into (\ref{t+1-t}), we obtain
\begin{equation}
    (t+1)^{\hat{n}}-t^{\hat{n}} = \sum_{s=1}^{\hat{n}} \genfrac\{\}{0pt}{}{\hat{n}}{s}s \sum_{j=0}^{t-1} \beta_j(t-j-1)_{s-1}.
\end{equation}
Thus, it follows from (\ref{strt-j}) that the left hand side of (\ref{t-jt+1}) equals to the right hand side. Hence, $f_0(t) \equiv 0$ and we can conclude that (\ref{uprule}) is equivalent to (\ref{upruletilx}).
\hfill $\blacksquare$

\begin{remark}
    It follows immediately from this lemma that an algorithm of the form (\ref{uprule}) satisfying Condition \ref{assumpbeta} will achieve convergence with zero steady-state error for all $f(x,t) \in \mathscr{F}^q_{m,L}$ if and only if the corresponding algorithm (\ref{upruletilx}) is such that $\tilde{x}(t)\rightarrow0$ as $t\rightarrow\infty$ for all $\tilde{f}(\tilde{x}) \in \mathscr{F}^q_{m,L}$.
\end{remark}

We can now reformulate (\ref{upruletilx}) as an uncertain linear feedback system:
\begin{align}
\label{sys}
   \chi(t+1) =& A\chi(t) + \tilde{B} \tilde{u}(t), \nonumber\\
    \tilde{y}(t) =& \tilde{C} \chi(t), \nonumber\\
    \tilde{u}(t) =& -\bold{\Delta} \tilde{y}(t), 
\end{align}
where $\tilde{\chi}(t) = [\tilde{x}(t-k)^T \quad \dots \quad \tilde{x}(t)^T]^T \in \mathbb{R}^{(k+1)p}$, $A \in \mathbb{R}^{(k+1)p\times (k+1)p}, \tilde{B} \in \mathbb{R}^{(k+1)p\times (k+1)p}, \tilde{C} \in \mathbb{R}^{(k+1)p\times (k+1)p}$,
\begin{equation}
\label{itermatr}
     A = A_0 \otimes I_p, \quad
     \tilde{B} = \tilde{B}_0 \otimes I_p, \quad  \tilde{C} = \tilde{C}_0 \otimes I_p,
\end{equation}
and
\begin{align}
\label{A0B0C0}
    A_0 &= {\left[ 
    \begin{array}{c|c}
         \textbf{0} & I_k  \nonumber\\
         \hline
         -\beta_{k-1}& \begin{matrix}
             a
         \end{matrix}
    \end{array}\right],}\nonumber\\
    \tilde{B}_0&= \left[ \begin{array}{c}
         \textbf{0}  \nonumber\\
         \hline
         \begin{matrix}
             \alpha_k &\alpha_{k-1}&\dots&\alpha_0
         \end{matrix} 
    \end{array} \right],\nonumber\\
    \tilde{C}_0&=  \begin{array}{c}
         I_{k+1} 
    \end{array}.
\end{align}
Here,
\begin{equation*}
    a = {\footnotesize\left[\begin{matrix}
              \beta_{k-1}-\beta_{k-2} & \beta_{k-2}-\beta_{k-3} &\dots&\beta_1-\beta_0 &1+\beta_0
          \end{matrix}\right]}
\end{equation*}
and $\boldsymbol{\Delta}$ is a block diagonal matrix with $k+1$ copies of the symmetric matrix $\Delta$ on its diagonal, 
\begin{equation}
\label{deltamatrix}
    \bold{\Delta} = \begin{bmatrix}
        \Delta & & \\ &\Delta \\& &\ddots & \\& & &\Delta
    \end{bmatrix},
\end{equation}
where $\Delta$ satisfies (\ref{classfunc}), $\otimes$ denotes the Kronecker product, $I_k$ is  the $k \times k$ identity matrix, and \textbf{0} denotes the zero matrix or vector of appropriate dimension.

Define $\hat{A}(\bold{\Delta})=A-\tilde{B} \bold{\Delta} \tilde{C}$, where $A, \tilde{B}, \tilde{C}$ are defined as in (\ref{itermatr}). We also define $\hat{A_0}(\lambda) = A_0- \lambda \tilde{B}_0 \tilde{C}_0$ where $m \leq \lambda \leq L$.

The following lemma, which follows from \cite[Theorem~10.1.4, p.~301]{IterativeSolution} and \cite[Theorem~1]{UGRINOVSKII2023111129}, relates the convergence rate of the algorithm (\ref{upruletilx}) to the spectral radius of $\hat{A_0}(\lambda)$.

\begin{lemma}
\label{convergespectral}
    Given any $\alpha$, $\beta$, $m$, $L$, such that the algorithm (\ref{upruletilx}) is globally convergent with the convergence rate
    \begin{equation}
        r(\alpha,\beta,m,L) \triangleq \sup_{mI \leq \Delta \leq LI} \rho(\hat{A}(\Delta)) \label{matrixr}.
    \end{equation}
    Then
    \begin{equation}
        r(\alpha,\beta,m,L)=\sup_{m \leq \lambda \leq L} \rho(\hat{A_0}(\lambda)). \label{scalarr}
    \end{equation}
    Here, $\rho(\cdot)$ denotes the spectral radius of a matrix.
\end{lemma}

 In order to apply this lemma to calculate the convergence rate of the algorithm (\ref{upruletilx}) and hence the algorithm (\ref{uprule}), we first consider the characteristic polynomial of $\hat{A_0}(\lambda)$. Indeed, we compute $\hat{A_0}(\lambda)$ as

    \begin{align}
        \hat{A_0}(\lambda) &= A_0- \lambda \tilde{B}_0 \tilde{C}_0 \nonumber \\
        &= \begin{bmatrix}
            0 & 1 & \dots & 0\\
            \vdots & \vdots & \ddots & \vdots\\
            0 & 0 & \dots & 1\\
            \eta_{k-1} & \eta_{k-2}& \dots &\eta_{0}
        \end{bmatrix}
    \end{align}
    where
    \begin{align*}
        \eta_{k-1} &= -\beta_{k-1}-\alpha_k \lambda,\\
        \eta_{k-2} &=\beta_{k-1}-\beta_{k-2}-\alpha_{k-1}\lambda,\\
        &\vdots\\
        \eta_0 &=1+\beta_0-\alpha_0 \lambda.
    \end{align*}
    Then we compute the characteristic equation of $\hat{A_0}(\lambda)$ as
    \begin{equation}
    \label{charA}
        (z-1) \left(z^k - \sum_{j=0}^{k-1} \beta_j z^{k-j-1} \right)+ \lambda \sum_{j=0}^{k}\alpha_j z^{k-j}=0.
    \end{equation}
    Equation (\ref{charA}) can be written as $1+\lambda \tilde{P}(z)\tilde{K}(z)=0$, where
    \begin{align}
        \tilde{P}(z)&\triangleq \frac{1}{z-1},\nonumber\\
        \tilde{K}(z)&\triangleq \frac{\tilde{N}(z)}{\tilde{D}(z)},\nonumber\\
        \tilde{N}(z)& \triangleq \sum_{j=0}^{k}\alpha_j z^{k-j},\nonumber\\
        \tilde{D}(z) &\triangleq z^k - \sum_{j=0}^{k-1} \beta_j z^{k-j-1}.\label{tilD}
    \end{align}

The following lemma relates Condition \ref{assumpbeta} to the condition that the algorithm (\ref{uprule}) incorporates at least $n$ discrete time integrations.
\begin{lemma}
    Condition \ref{assumpbeta} is satisfied if and only if the polynomial $\tilde{D}(z)$ defined in (\ref{tilD}) has at least $n-1$ roots at $z=1$.
\end{lemma}

\textit{Proof:}
Let $\hat{k}\geq k$ be given and let $\bar{k} = \hat{k}-k \geq 0$. Also, define
\begin{align*}
    \hat{D}(z) &= z^{\hat{k}} - \sum_{j=0}^{k-1} \beta_j z^{\hat{k}-j-1}\\
    &=z^{k}z^{\bar{k}} -  \sum_{j=0}^{k-1} \beta_j z^{k+\bar{k}-j-1}\\
    &=z^{k}z^{\bar{k}} -  \sum_{j=0}^{k-1} \beta_j z^{\bar{k}} z^{k-j-1}\\
    &=z^{\bar{k}} \tilde{D}(z).
\end{align*}
From this, it follow that $\tilde{D}(z)$ will have at least $n-1$ roots at $z=1$ if and only if $\hat{D}(z)$ has at least $n-1$ roots at $z=1$.

We now consider the first $n-1$ derivatives of $\hat{D}(z)$:
\begin{equation}
    \hat{D}^{(r)}(z) = (\hat{k})_r z^{\hat{k}-r} - \sum_{j=0}^{k-1}\beta_j(\hat{k}-j-1)_r z^{\hat{k}-j-1-r},
\end{equation}
for $r = 0,1,\dots,n-1$.
Then $\hat{D}(z)$ will have at least $n-1$ roots at $z=1$, if and only if
    \begin{align*}
        \hat{D}^{(r)}(1) &= (\hat{k})_r - \sum_{j=0}^{k-1}\beta_j(\hat{k}-j-1)_r \\
        &=0, 
    \end{align*}
for $r = 0,1,\dots,n-1$.
This is equivalent to
\begin{equation}
    (\hat{k})_r =  \sum_{j=0}^{k-1}\beta_j(\hat{k}-j-1)_r,
\end{equation}
for $r = 0,1,\dots,n-1$.
This is Condition \ref{assumpbeta}, thereby verifying the lemma. \hfill $\blacksquare$

\section{Main Results} \label{sec3}
The main result of this paper, which is a fundamental convergence rate bound for algorithms of the form (\ref{uprule}), is stated in the following theorem. 

\begin{theorem}
\label{maintheo1:r>=phb}
    Let $L>m>0$ and $\kappa = \frac{L}{m}$ be given. Then any algorithm $\Sigma(\alpha,\beta,m,L)$ of the form (\ref{uprule}) which is globally convergent with zero steady-state error for any $f(x,t) \in \mathscr{F}_{m,L}^q$, will have a convergence rate $r(\alpha,\beta,m,L)$ satisfies:
    \begin{align}
    \label{maxrateineq}
        r(\alpha,\beta,m,L) \geq \sqrt[n]{\rho_{\text{\scriptsize HB}}}
        &= \left(\frac{\sqrt{\kappa} - 1 }{\sqrt{\kappa} + 1}\right)^{\frac{1}{n}} \nonumber \\ \quad &= \left(\frac{\sqrt{L} - \sqrt{m}}{\sqrt{L} + \sqrt{m}}\right)^{\frac{1}{n}}.
    \end{align}
\end{theorem}
The following corollary follows immediately from the above theorem.
\begin{corollary}
  Let $L>m>0$ and $\kappa = \frac{L}{m}$ be given. Then any algorithm $\Sigma(\alpha,\beta,m,L)$ which is globally convergent with zero steady-state error for any set of cost functions $\mathscr{F} \supset \mathscr{F}_{m,L}^q$, will have a convergence rate which satisfies the bound (\ref{maxrateineq}).    
\end{corollary}

The following theorem shows that the bound given in Theorem \ref{maintheo1:r>=phb} is tight and gives the construction of an algorithm $\Sigma(\alpha^*,\beta^*,m,L)$ which achieves the bound.
\begin{theorem}
\label{maintheo}
Let $L>m>0$ and $\kappa = \frac{L}{m}$ be given. Then there exists an algorithm $\Sigma(\alpha^*,\beta^*,m,L)$ of the form (\ref{uprule}) which is globally convergent with zero steady-state error for $f(x,t) \in \mathscr{F}^q_{m,L}$ and has the convergence rate 
\begin{align}
\label{maintheorate}
    r(\alpha^*,\beta^*,m,L)&=\sqrt[n]{\rho_{\text{\scriptsize HB}}} \nonumber\\
    &=\left(\frac{\sqrt{\kappa} - 1 }{\sqrt{\kappa} + 1}\right)^{\frac{1}{n}}= \left(\frac{\sqrt{L} - \sqrt{m}}{\sqrt{L} + \sqrt{m}}\right)^{\frac{1}{n}}.
\end{align}
In this algorithm, $k=2n-1$ and the corresponding values of $\alpha^*$ and $\beta^*$ are given by
    \begin{align}
    \label{eq:alpha}
        \alpha_j 
        &= \tilde K_1 \binom{2n}{j}(-\rho^2)^j
          + \tilde K_2 \binom{2n}{j}(-1)^j \nonumber
          \\& \quad- \tilde K_3
            \sum_{i=0}^{j}
              \binom{n}{i}\,\binom{n}{j-i}\,(-\rho^2)^i\,(-1)^{\,j-i},
    \end{align}
for $j = 0,1,\dots,2n,$ and 

    \begin{equation}
    \label{eq:beta}
    \beta_j
    = -\sum_{i=0}^{j+1}
         \binom{n}{i}\,\binom{n-1}{\,j+1-i\,}\,
         (-\rho^2)^i\,(-1)^{\,j+1-i},
    \end{equation}
for $j = 0,1,\dots,2n-1,$ where
    \begin{align}
    \label{correstilK}
    \tilde{K}_1&=\frac{(\sqrt{L}+\sqrt{m})^2}{4Lm},\nonumber\\
    \tilde{K}_2&=\frac{(\sqrt{L}-\sqrt{m})^2}{4Lm},\nonumber\\
    \tilde{K}_3&=\frac{L+m}{2Lm}.
    \end{align}
\end{theorem}
\smallskip

In order to prove the above theorems, we introduce a series of lemmas.

\begin{lemma}
\label{ztransx*(t)}
    The signal $x^*(t)$ in (\ref{x*t}) has a $z$-transform of the form
    \begin{equation}
    \label{ztransx*}
        X^*(z) = a_0\frac{z}{z-1}+a_1\frac{B_1(z)}{(z-1)^2}+\dots + a_{n-1}\frac{B_{n-1}(z)}{(z-1)^n},
    \end{equation}
    where $B_r(z)$ is a polynomial of order $r+1$ such that $B_r(1)=r!$ for $r = 0,1,\dots,n-1$.
\end{lemma}

\textit{Proof:}
    The proof proceeds by induction on $r$.
    By inspection, 
    \begin{equation*}
        \mathscr{Z}(t^0) = \frac{z}{z-1} = \frac{B_0(z)}{z-1},
    \end{equation*}
    where $B_0(z) = z$ is a polynomial of order 1 and $B_0(1) = 1$. Hence, the base case holds.
    
    Now assume that for a given $r \geq 0$, the following formula holds
    \begin{equation}
    \label{Zt^n}
        \mathscr{Z}(t^r) = \frac{B_r(z)}{(z-1)^{r+1}}
    \end{equation}
    where $B_r(z)$ is a polynomial of order $r+1$ and $B_r(1)=r!$. 
   We now consider the case of $t^{r+1}$ and use a property of $z$-transforms; e.g., see \cite[Equation 2.9, p.~16]{fadali2009digital}:
    \begin{align}
        \mathscr{Z}(t^{r+1}) &= \mathscr{Z}(t t^r) \nonumber\\
        &=-z\frac{d}{dz}\mathscr{Z}(t^r) \nonumber\\
        &=-z\frac{d}{dz}\frac{B_r(z)}{(z-1)^{r+1}} \nonumber\\
        &=-z \left( \frac{(z-1)^{r+1}B'_r(z)-(r+1)(z-1)^rB_r(z)}{(z-1)^{2r+2}} \right) \nonumber\\
        &= \frac{z(r+1)B_r(z)-z(z-1)B'_r(z)}{(z-1)^{r+2}} \nonumber\\
        &=\frac{B_{r+1}(z)}{(z-1)^{r+2}}, \label{Zt^n+1}
    \end{align}
    where $B'_r(z)$ denotes the derivative of $B_r(z)$.
    Using (\ref{Zt^n+1}), and the fact that $B_r(z)$ is a polynomial of order $r+1$, it follows that $B_{r+1}(z)$ is a polynomial of order $r+2$ and
    \begin{equation}
        B_{r+1}(1) =(r+1)B_r(1) = (r+1)r! = (r+1)!.
    \end{equation}
    Thus by induction, (\ref{Zt^n}) holds for all $r\geq0$. From this, the lemma follows.
\hfill $\blacksquare$

The following lemma provides a necessary condition for achieving zero steady-state error, and amounts to a version of the internal model principle for the problem under consideration.
\begin{lemma}
    \label{maintheonintegrators}
    Given any $L>m>0$, suppose an algorithm $\Sigma(\alpha,\beta,m,L)$ of the form (\ref{uprule}) is globally convergent with zero steady-state error for any $f(x,t) \in \mathscr{F}_{m,L}^q$. Then 
    \begin{align}
        D(z) &\triangleq (z-1)\tilde{D}(z)\nonumber\\ &=(z-1)\left( z^k - \sum_{j=0}^{k-1} \beta_j z^{k-j-1} \right)
    \end{align}
    has at least $n$ roots at $z=1$. That is, $\tilde{D}(z)$ has at least $n-1$ roots at $z=1$.
\end{lemma}

\textit{Proof:}
In order to prove this lemma, we write the iterative algorithm (\ref{uprule}) in terms of a Luré system by letting 
\begin{equation*}
    \chi(t) = \begin{bmatrix}
        x(t-k)^T& \dots& x(t)^T
    \end{bmatrix}^T \in \mathbb{R}^{(k+1)p}
\end{equation*}
be the state vector, and defining the reference signal
\begin{equation*}
    \chi^*(t) = \begin{bmatrix}
        x^*(t-k)^T&\dots&x^*(t)^T
    \end{bmatrix}^T \in \mathbb{R}^{(k+1)p}.
\end{equation*}
Also, let $A$, $B$ and $C$ be matrices defined by 
\begin{equation}
    \label{newsysmatrix}
    A = A_0 \otimes I_p, \quad B = B_0 \otimes I_p, \quad  C = C_0 \otimes I_p,
\end{equation}
with $A_0 \in \mathbb{R}^{(k+1) \times (k+1)}$ defined in (\ref{A0B0C0}),
 \begin{align}
    B_0&= \left[ \begin{array}{c}
         \textbf{0}  \\
         \hline
            1
    \end{array} \right] \in \mathbb{R}^{(k+1)}, \nonumber\\
    C_0&= \left[ \begin{array}{c}
        \begin{matrix}
            \alpha_k &\alpha_{k-1}&\dots&\alpha_0
        \end{matrix}
    \end{array}
      \right]\in \mathbb{R}^{(k+1)}. \label{b0c0notil}
\end{align}
It follows that for any $f(x,t)\in\mathscr{F}_{m,L}^q$, (\ref{uprule}) can be rewritten as
\begin{align}
    \label{Xt+1orig}
    \chi (t+1) &= A\chi(t) - B\diag{(\bold{\Delta)}} C (\chi(t)-\chi^*(t)) 
\end{align}
where  $\bold{\Delta}$ is defined as in (\ref{deltamatrix}).
Since we have assumed that the algorithm guarantees convergence with zero steady-state error for any $f(x,t)\in\mathscr{F}_{m,L}^q$, it follows that we will obtain convergence with zero steady-state error for $ f(x,t)\in\mathscr{F}_{m,L}^q$
with $\Delta = \lambda I_p$ and $m \leq \lambda \leq L$. In this case, $\bold{\Delta} = \lambda I$ and we can rewrite (\ref{Xt+1orig}) as
\begin{equation}
    \chi(t+1) = A\chi(t) - \lambda B C (\chi(t)-\chi^*(t)).
\end{equation}
This system, can be rewritten in state space form as
\begin{align}
    \label{newsys}
    \chi(t+1) &= A\chi(t) + Bu(t), \nonumber\\
    y(t) &= C\chi(t), \nonumber\\
    \xi(t) &= C\chi^*(t),\nonumber\\
    u(t) &= \lambda(\xi(t)-y(t)).
\end{align}

The system (\ref{newsys}) is illustrated in Figure \ref{fig:bdimp}.
\tikzstyle{block} = [draw, rectangle, 
    minimum height=3em, minimum width=6em]
\tikzstyle{sum} = [draw, circle, node distance=1cm]
\tikzstyle{input} = [coordinate]
\tikzstyle{output} = [coordinate]
\tikzstyle{pinstyle} = [pin edge={to-,thin,black}]

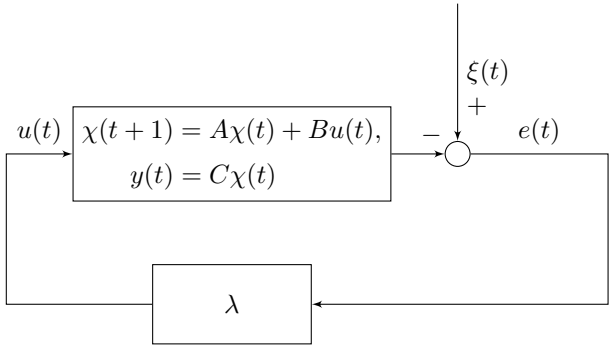
\begin{figure}[htp!]
    \centering
    \begin{tikzpicture}[auto, node distance=2cm,>=latex']
    \node [input, name=input] {};
    \node [block, right of=input,node distance=3cm] (controller) {$\begin{aligned}
       \chi(t+1) &= A \chi(t)+Bu(t),\\[2pt]
      y(t)      &= C \chi(t)
    \end{aligned}$};
    \node [sum, right of=controller, node distance=3cm] (disturbance) {};
    \node [output, right of=disturbance] (output) {};
    \draw [->] (controller) -- node [near end] {$-$} (disturbance);
    \node[block,below of=controller](u){$\lambda$};
    \draw [->] (disturbance) -- node{$e(t)$}(output) |- (u);
    \node[input, above of=disturbance](dist_input){$\xi(t)$};
    \draw[->]  (dist_input) --node{$\xi(t)$} node [near end] {$+$}(disturbance);
    \draw[->] (u) -| (input) -- node{$u(t)$}(controller);
\end{tikzpicture}
    \caption{Block diagram representation of the system (\ref{newsysmatrix}), (\ref{b0c0notil}), (\ref{newsys}).}
    \label{fig:bdimp}
\end{figure}

Let \(G(z)\) denote the transfer function $G(z) = C(zI-A)^{-1}B$, which takes the form
\[
G(z) = G_0(z)\otimes I_p,
\]
where the scalar transfer function $G_0(z)$ is given by
\begin{equation}
\label{newtranferfunc}
    G_0(z)=\frac{\alpha_0z^{k}+\alpha_1z^{k-1}+\dots+\alpha_{k-1}z+\alpha_k}{D(z)},
\end{equation}

where the denominator polynomial is
\begin{align*}
    D(z)&=z^k+(-1-\beta_0)z^{k-1}+(\beta_0-\beta_1)z^{k-2}\\
    &\quad+\dots+(\beta_{k-3}-\beta_{k-2})z^2+(\beta_{k-2}-\beta_{k-1})z\\&\quad+\beta_{k-1} \\&=(z-1)\left( z^k - \sum_{j=0}^{k-1} \beta_j z^{k-j-1} \right)\\&=(z-1) \tilde{D}(z).
\end{align*}
Note that (\ref{newtranferfunc}) follows from the fact that the system (\ref{newsysmatrix}), (\ref{b0c0notil}), (\ref{newsys}) is in controllable canonical form.

It follows from Figure \ref{fig:bdimp} that the error signal $e(t)$ is such that
\[
E(z) = -\lambda E(z)G(z)+\Xi(z),
\]
leading to
\begin{align*}
    E(z) &= (I+\lambda G(z))^{-1}\Xi(z)\\
    &=(I_p+\lambda G_0(z)\otimes I_p)^{-1}\Xi(z)
\end{align*}

where \(\Xi(z)\) is the \(z\)-transform of \(\xi(t)\) and $E(z)$ is the $z$-transform of $e(t)$.

Thus,
\[
E(z) = \frac{\Xi(z)}{1+\lambda G_0(z)},
\]
with
\[
\xi(t)=\alpha_k x^*(t-k)+\alpha_{k-1} x^*(t-k+1)+\dots+\alpha_0x^*(t).
\]

Applying the \(z\)-transform to \(\xi(t)\), we obtain
\[
\Xi(z)=\left(\alpha_k z^{-k}+\alpha_{k-1} z^{-k+1}+\dots+\alpha_0\right)x^*(z).
\]

It follows from Lemma \ref{ztransx*(t)} that 
\begin{align*}
    \Xi(z)&=\left(\alpha_k z^{-k}+\alpha_{k-1} z^{-k+1}+\dots+\alpha_0\right) \nonumber\\ &\quad \times\left( a_0\frac{B_0(z)}{z-1}+a_1\frac{B_1(z)}{(z-1)^2}+\dots+a_{n-1}\frac{B_{n-1}(z)}{(z-1)^n}\right).
\end{align*}
Hence,
\begin{align}
    \label{E(z)1}
    &E(z) = \frac{\alpha_k z^{-k}+\alpha_{k-1} z^{-k+1}+\dots+\alpha_0}{1+\lambda G_0(z)} \nonumber\\
    &\quad \times\left(a_0\frac{B_0(z)}{z-1}+a_1\frac{B_1(z)}{(z-1)^2}+\dots+a_{n-1}\frac{B_{n-1}(z)}{(z-1)^n}\right).
\end{align}

Multiplying the numerator and denominator of this expression by $(z-1)^n$, we rewrite (\ref{E(z)1}) as
\begin{align}
\label{Ez}
    E(z)&=\frac{\left(\alpha_k z^{-k}+\dots+\alpha_0\right)}{(z-1)^n+\lambda(z-1)^n G_0(z)} \nonumber \\
    &\quad \times \left((z-1)^{n-1}a_0B_0(z)+\dots+a_{n-1}B_{n-1}(z)\right)
\end{align}

Using the Final Value Theorem \cite[Theorem 2.1]{fadali2009digital}, the steady-state error is
\begin{equation*}
    \lim_{t\to\infty} e(t)=\lim_{z\to 1}(z-1) E(z).
\end{equation*}

Therefore, to achieve zero steady-state error, it is required that
\begin{equation*}
    \lim_{z\to 1}(z-1)E(z)=0.
\end{equation*}
Using (\ref{Ez}), this implies that the denominator of $G_0(z)$ must contain at least $n$ poles at $z=1$. That is, $D(z)$ has at least $n$ roots at $z=1$. Hence, $\tilde{D}(z)$ has at least $n-1$ roots at $z=1$.\hfill $\blacksquare$

The proofs of our main results use results on the optimal gain margin problem as established in \cite[Problem~2.4]{Tannenbaum1985} and incorporate a related special case from the theoretical development in \cite{UGRINOVSKII2023111129}. Generalizing the case in \cite{UGRINOVSKII2023111129}, our approach considers a plant that includes a discrete time $n$-th order integrator, enabling the solution of the polynomially time-varying optimization problem with zero steady-state error. The required result is contained in the following lemma.
\begin{figure}[htp!]
    \centering
    \begin{tikzpicture}[auto, node distance=2cm,>=latex']
    \node [input, name=input]{$U(z)$};
    \node [sum, right of=input,node distance=1.2cm] (sum) {};
    \node [block, right of=sum,node distance=1.7cm] (compensator) {$K(z)$};
    \node [block, right of=compensator, node distance=2.7cm] (plant) {$\lambda P(z)$};
    \node [output, right of=disturbance,node distance=2cm] (output) {};
    \coordinate [below of=compensator] (measurements) {};
    
    \draw [->] (input) -- node{$U(z)$} node[pos=0.95]{$+$}(sum);
    \draw [->] (sum) -- node{}(compensator);
    \draw [->] (compensator) -- (plant);
    \draw [->] (plant) -- node [name=y]{$Y(z)$} (output);
    \draw [-] (y) |- (measurements);
    \draw [->] (measurements) -| node[pos=0.99] {$-$} (sum);
    \end{tikzpicture}
    \caption{Linear uncertain system.}
    \label{fig:BD}
\end{figure}
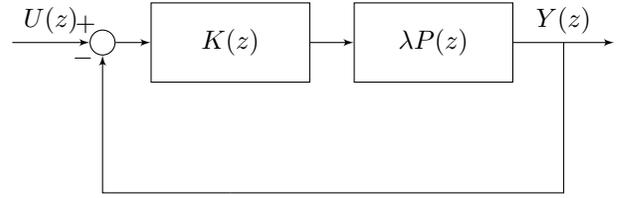
\begin{lemma}
\label{lemma1}
    Consider the uncertain linear feedback system depicted in Figure \ref{fig:BD}, which comprises the plant $P(z) = \frac{1}{(z-1)^n}$, which is a discrete time $n$-th order integrator, an uncertain constant gain $\lambda$ and a compensator $K(z)$. Let $\rho \in (0,1)$ be given and suppose $P(z)K(z)$ is strictly proper. If the compensator $K(z)$ places all of the poles of this system in the disk $|z|<\rho$ for all $\lambda \in [m,L]$ then 
    \begin{equation}
    \label{rhorhohb}
        \rho \geq \sqrt[n]{\rho_{\text{\scriptsize HB}}},
    \end{equation}
    where $\rho_{\text{\scriptsize HB}}$ is the convergence rate of the heavy ball method (\ref{rateHB}).
\end{lemma}

\textit{Proof:}
    The proof of the lemma follows along the same lines as the proof of \cite[Lemma~2]{UGRINOVSKII2023111129}, with modifications to accommodate our problem setting.
    
    We introduce the sensitivity function $S(z)$ defined as
    \begin{equation}
    \label{sensifunc}
        S(z) = \left( 1+ \frac{m+L}{2} P(z)K(z) \right)^{-1}.
    \end{equation}
    The $n$ poles of $P(z)$ at $z=1$, are $n$ zeros of $S(z)$. Also the zeros of $P(z)$ at $z=\infty$, correspond to a zero of $1-S(z)$, since $P(z)K(z)$ is strictly proper. Therefore, we conclude that
        \begin{equation}
        \label{interpoint1}
            S(1) = 0,
        \end{equation}
    with multiplicity of $n$,
    and
    \begin{equation}
    \label{interpoint2}
        S(\infty)=1.
    \end{equation}

    Define $\mathscr{G} \triangleq \left( -\infty,\frac{2m}{m-L} \right] \cup \left(\frac{2L}{L-m}, +\infty \right)$ and $\mathscr{G}^\textbf{C} = \mathbb{C} \setminus \mathscr{G}$, where $\mathbb{C}$ is the set of complex numbers.
    
\begin{figure}[htb!]
        \centering
        \includegraphics[width=0.5\linewidth]{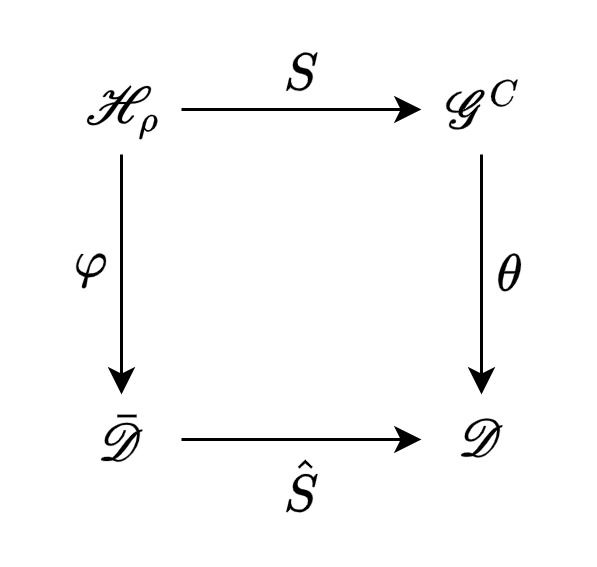}
        \vspace{-0.5cm}
        \caption{Commutative diagram}
        \label{fig:map}
    \end{figure}

    It follows from \cite[Lemma~2]{UGRINOVSKII2023111129} and \cite[Lemma~2.3]{Tannenbaum1985}, that in order to place the poles of the closed loop system within the interior of the open disk $|z|<\rho$ for all $\lambda \in [m,L]$, it is equivalent to demonstrate the existence of a sensitivity function $S(z)$ defined in (\ref{sensifunc}) which is analytic in $\mathscr{H}_\rho \triangleq \{ |z| \geq \rho\} \cup \{ \infty \}$ that satisfies the conditions (\ref{interpoint1}), (\ref{interpoint2}) and $S(z)\in \mathscr{G}^\textbf{C}$ for all $z\in \mathscr{H}_\rho$. To determine the existence of such an $S(z)$, we reformulate the problem as a Nevanlinna Pick interpolation problem, as described in \cite{UGRINOVSKII2023111129,Tannenbaum1985}, using the commutative diagram in Figure \ref{fig:map}.
    
    In this figure,
    \begin{align*}
        \mathscr{\Bar{\mathscr{D}}} &= \{z:|z| \leq 1\}, \\
        \mathscr{\mathscr{D}} &= \{z:|z| < 1\},
    \end{align*}
    and the function $\varphi(z) = \rho z^{-1}$ maps $\mathscr{H}_\rho$ into $\bar{\mathscr{D}}$. Similarly, it follows from \cite{UGRINOVSKII2023111129,Tannenbaum1985} that the function 
    \begin{equation}
        \theta(z) = \frac{1-\sqrt{\dfrac{1-\frac{z (L-m)}{2 L}}{1-\frac{z (m-L)}{2 m}}}}{1+\sqrt{\dfrac{1-\frac{z (L-m)}{2 L}}{1-\frac{z (m-L)}{2 m}}}}
    \end{equation}
    is analytic in $\mathscr{G}^\textbf{C}$, maps $\mathscr{G}^\textbf{C}$ into $\mathscr{D}$. Also, the function 
    \begin{equation}
    \label{hatS=theta}
        \hat{S}(z) = \theta(S(\varphi(z)))
    \end{equation}
    should be constructed to be analytic in $\mathscr{\Bar{\mathscr{D}}}$, and map the closed unit disk $\mathscr{\Bar{\mathscr{D}}}$ into the open unit disk $\mathscr{\mathscr{D}}$.
    
    As illustrated in the commutative diagram in Figure \ref{fig:map}, the existence of a function $S(z)$ satisfying (\ref{interpoint1}) and (\ref{interpoint2}) which is analytic in $\mathscr{H}_\rho$, and maps $\mathscr{H}_\rho$ into $\mathscr{G}^\textbf{C}$, is equivalent to a Nevanlinna Pick interpolation problem for the function $\hat{S}(z)$.
    
    The Nevanlinna Pick interpolation problem is defined as follows. Given a set of points $\{ z_i \}$ in the closed unit disk $\mathscr{\Bar{\mathscr{D}}}$ and a corresponding set of target values $\{w_i \}$, find a function $\hat{S}(z)$ which is analytic in $\mathscr{D}$ such that $|\hat{S}(z)| < 1 $ for all $z\in  \mathscr{\mathscr{D}}$, where $\mathscr{\mathscr{D}}$ is the open unit disk, and satisfying the interpolation conditions $\hat{S}(z_i)=w_i$ with specified multiplicities.  
    
    It follows from (\ref{interpoint1}) and (\ref{interpoint2}) that the required interpolation values for $\hat{S}(z)$ are
    \begin{equation}
        \hat{S}(\rho) =\theta(S(1))=0\label{hatsrhointerpoint} 
    \end{equation}
    with multiplicity $n$, and
    \begin{equation}
    \label{hats0interpoint}
        \hat{S}(0) = \theta(S(\infty))=\theta(1).
    \end{equation}
    Note that it is straightforward to verify that $\theta(1) = \rho_{\text{\scriptsize HB}}$.
    
    Handling an interpolation condition with multiplicities requires imposing conditions on the derivatives of $\hat{S}(z)$, which complicates the interpolation problem. To avoid this complexity, we introduce small perturbations to the plant's poles by redefining the plant transfer function as
    \begin{equation}
        P(z) = \frac{1}{\prod_{i=1}^{n}(z-1-\varepsilon_i)}
    \end{equation}
     where $n$ denotes the number of integrators with $n>1$ and the \( \varepsilon_i \) are small positive parameters. As $\varepsilon_i \rightarrow 0$, the plant $P(z)$ reduces to $\frac{1}{(z-1)^n}$, which represents a $n$-th order integrator. This modification allows us to avoid multiplicities in the interpolation conditions by ensuring that the poles and zeros of $P(z)$ are distinct, thereby simplifying the Nevanlinna-Pick interpolation problem. In this modified problem, the poles of $P(z)$ at $z=1+\varepsilon_i$, are zeros of $S(z)$. Therefore, the interpolation conditions of $S(z)$ become
        \begin{align}
        &S(1+\varepsilon_i)=0, \quad i = 1,2,\dots,n,\\
        &S(\infty)=1, \quad
        \end{align}
    all with multiplicity one.
    By perturbing the interpolation points in this manner, we obtain modified interpolation conditions on $\hat{S}(z)$, which yields
        \begin{align}
        &\hat{S}\left(\frac{\rho}{1+\varepsilon_i}\right) = \theta(S(1))=\theta(0) =0,\label{modifyinterpoint1}\\
        &\hat{S}(0) = \theta(S(\infty))=\theta(1),\label{modifyinterpoint2}
        \end{align}
    all with multiplicity one.
    
    Now we are in the position to determine the existence of such an $\hat{S}(z)$. By \cite[Theorem~2.2]{BoundAnalyticFunc} there exists an interpolation function $\hat{S}(z)$ with $\hat{S}(z) <1$ for all $z\in \mathscr{D}$ that satisfies (\ref{modifyinterpoint1})-(\ref{modifyinterpoint2}) if and only the Pick matrix $M\in \mathbb{R}^{n+1\times n+1}$ written as
    \begin{equation}
    \label{pickmatr}
        M = \left[\begin{array}{c|c}
            P & Q \\
            \hline
            Q^T & R
        \end{array}\right]
    \end{equation}
    is positive semi-definite.
    In (\ref{pickmatr}), $P \in \mathbb{R}^{n \times n}$, 
    \begin{equation}
        P_{i,j} = \frac{1}{1-\frac{\rho^2}{(1+\varepsilon_{i})(1+\varepsilon_{j})}}, \quad i,j=1,...,n,
    \end{equation}
    $Q=\textbf{1}$ is a column vector of ones and $R$ is a scalar, $R=1-\theta(1)^2$.
    To satisfy the positive semi-definiteness of the Pick matrix $M$, a necessary condition is that its determinant be non-negative. To compute the determinant of the Pick matrix $M$, we employ the Schur complement \cite{Schur1917}, which yields
    \begin{equation}
    \label{detM}
        \det(M) = \det(R) \det(P-QR^{-1}Q^{T}).
    \end{equation}
    Since $R$ is a scalar
    \begin{equation}
    \label{R-1}
        R^{-1} = \frac{1}{1-\theta(1)^2},
    \end{equation}
    and
    \begin{equation}
        \label{detR}
        \det (R) = 1-\theta(1)^2.
    \end{equation}
    Let $T = P-QR^{-1}Q^{T}$ where $T \in  \mathbb{R}^{n \times n}$. Since $Q$ is a vector of ones, we obtain
    \begin{align*}    
        T_{i,j} &=  \frac{1}{1-\frac{\rho^2}{(1+\varepsilon_{i})(1+\varepsilon_{j})}}-\frac{1}{1-\theta(1)^2}\\
        &=\frac{1+\varepsilon_{i}}{(1+\varepsilon_{i})-\frac{\rho^2}{1+\varepsilon_{j}}}-\frac{1}{1-\theta(1)^2}\\
        &=\frac{(1+\varepsilon_{i})(1-\theta(1)^2)-((1+\varepsilon_{i})-\frac{\rho^2}{1+\varepsilon_{j}})}{(1+\varepsilon_{i})-\frac{\rho^2}{1+\varepsilon_{j}}} \times \frac{1}{1-\theta(1)^2}\\
        &=\frac{(1+\varepsilon_{i})(1-\theta(1)^2)-(1+\varepsilon_{i})+\frac{\rho^2}{1+\varepsilon_{j}}}{(1+\varepsilon_{i})-\frac{\rho^2}{1+\varepsilon_{j}}} \times \frac{1}{1-\theta(1)^2}\\
        &=\frac{-\theta(1)^2(1+\varepsilon_{i})+\frac{\rho^2}{1+\varepsilon_{j}}}{(1+\varepsilon_{i})-\frac{\rho^2}{1+\varepsilon_{j}}} \times \frac{1}{1-\theta(1)^2}.
    \end{align*}
    Let $T = UV$, where $U \in  \mathbb{R}^{n \times n}$,
    \begin{equation}
    \label{U}
        U_{i,j} = \frac{-\theta(1)^2(1+\varepsilon_{i})+\frac{\rho^2}{1+\varepsilon_{j}}}{(1+\varepsilon_{i})-\frac{\rho^2}{1+\varepsilon_{j}}}
    \end{equation}
    and $V \in  \mathbb{R}^{n \times n}$ is a diagonal matrix whose diagonal entries are given by $\frac{1}{1-\theta(1)^2}$.
    It is clear that 
    \begin{equation}
    \label{detT}
        \det(T) = \det(U)\det(V),
    \end{equation}
    and 
    \begin{equation}
    \label{detV}
        \det(V) = \left(\frac{1}{1-\theta(1)^2} \right)^n.
    \end{equation}
    Also, according to \cite[Theorem 1]{gendet},
    \begin{equation}
    \label{detU}
        \det(U) = (-1+\theta(1)^2)^{n-1} \frac{U_1 U_2}{\prod_{i,j}^n(1+\varepsilon_i-\frac{\rho^2}{1+\varepsilon_j})}
    \end{equation}
    where 
    \begin{align}
    \label{U1}
        U_1 &= -\prod_j^n\frac{-\rho^2}{1+\varepsilon_j}+(-1)^{n-1}(-\theta(1)^2)\prod_i^n(1+\varepsilon_i) \nonumber \\
        &=(-1)^{n+1}\prod_j^n\frac{\rho^2}{1+\varepsilon_j}+(-1)^{n}\theta(1)^2\prod_i^n(1+\varepsilon_i)
     \end{align}
     and
     \begin{align}
     \label{U2}
        U_2 &= \prod_{i<j}(1+\varepsilon_i-1-\varepsilon_j)(\frac{-\rho^2}{1+\varepsilon_i}-\frac{-\rho^2}{1+\varepsilon_j}) \nonumber\\
        &=\prod_{i<j}\frac{\rho^2(\varepsilon_i-\varepsilon_j)^2}{(1+\varepsilon_i)(1+\varepsilon_j)}.
    \end{align}
    Substituting (\ref{detV}), (\ref{detU}), (\ref{detT}) and (\ref{detR}) into (\ref{detM}) yields
    \begin{align}
    \label{detM1}
        \det(M) &= (1-\theta(1)^2)\left(\frac{1}{1-\theta(1)^2} \right)^n (-1+\theta(1)^2)^{n-1} \nonumber\\ &\quad \times \frac{U_1 U_2}{\prod_{i,j}^n(1+\varepsilon_i-\frac{\rho^2}{1+\varepsilon_j})} \nonumber \\
        &= (-1)^{n-1} \frac{U_1 U_2}{\prod_{i,j}^n(1+\varepsilon_i-\frac{\rho^2}{1+\varepsilon_j})}.
    \end{align}
    Substituting (\ref{U1}) into (\ref{detM1}) yields
    \begin{align}
        &\det(M) \nonumber\\  & = (-1)^{n-1} \nonumber\\ & \quad \times \frac{\left((-1)^{n+1}\prod_j^n\frac{\rho^2}{1+\varepsilon_j}+(-1)^{n}\theta(1)^2\prod_i^n(1+\varepsilon_i) \right)U_2}{\prod_{i,j}^n(1+\varepsilon_i-\frac{\rho^2}{1+\varepsilon_j})} \nonumber \\ \quad&=\frac{\left((-1)^{2n}\prod_j^n\frac{\rho^2}{1+\varepsilon_j}+(-1)^{2n-1}\theta(1)^2\prod_i^n(1+\varepsilon_i) \right)U_2}{\prod_{i,j}^n(1+\varepsilon_i-\frac{\rho^2}{1+\varepsilon_j})} \nonumber \\
        &=\frac{\left(\prod_j^n\frac{\rho^2}{1+\varepsilon_j}-\theta(1)^2\prod_i^n(1+\varepsilon_i) \right)U_2}{\prod_{i,j}^n(1+\varepsilon_i-\frac{\rho^2}{1+\varepsilon_j})}.
    \end{align}
    Since $\rho \in (0,1)$ and $\theta(1)\in (0,1)$, $P>0$ and $Q>0$ for all $\varepsilon_i>0$.
    Hence, a necessary condition for the Pick matrix $M$ to be positive semi-definite is $\det(M)\geq 0$. Letting $\varepsilon_i,\varepsilon_j \rightarrow 0$ from the positive direction, we obtain
    \begin{align}
        \lim_{\varepsilon_i, \varepsilon_j \rightarrow0} \frac{\det(M)}{U_2} &= \lim_{\varepsilon_i, \varepsilon_j \rightarrow 0} \frac{\left(\prod_j^n\frac{\rho^2}{1+\varepsilon_j}-\theta(1)^2\prod_i^n(1+\varepsilon_i) \right)}{\prod_{i,j}^n(1+\varepsilon_i-\frac{\rho^2}{1+\varepsilon_j})} \nonumber\\
        &= \frac{\rho^{2n}-\theta(1)^2}{
        (1-\rho^2)^{n^2}
    }.
    \end{align}
    It follows that a necessary condition for the Pick matrix $M$ to be positive semi-definite is 
    \begin{equation}
    \label{rho2theta}
        \frac{\rho^{2n}-\theta(1)^2}{
        (1-\rho^2)^{n^2}} \geq 0,
    \end{equation}
    which is equivalent to
    \begin{equation}
    \label{rhontheta}
        \rho^n \geq \theta(1).
    \end{equation}
    Since $\theta(1) = \rho_{\text{\scriptsize HB}}$, the condition (\ref{rhontheta}) is equivalent to (\ref{rhorhohb}). Thus, if there exists a compensator $K(z)$ satisfying the conditions of the lemma, then (\ref{rhorhohb}) must be satisfied.
\hfill $\blacksquare$

Now, we can prove Theorem \ref{maintheo1:r>=phb}.

\textit{Proof of Theorem \ref{maintheo1:r>=phb}:}
    Using Lemma \ref{maintheonintegrators}, a necessary condition for the algorithm $\Sigma(\alpha,\beta,m,L)$ to track a polynomially varying optimal point with zero steady-state error, for any $f(x,t) \in \mathscr{F}^q_{m,L}$, is that $\tilde{D}(z)$ defined in (\ref{tilD}) has at least $n-1$ roots at $z=1$. Moreover, Lemmas \ref{convergespectral} and \ref{lemma1} imply that the corresponding system in (\ref{sys}) must position all poles inside the disk defined by $|z|\leq \rho$ for every $\lambda \in [m,L]$ in order to achieve a convergence rate $\rho$. In order to meet this requirement, condition (\ref{rhorhohb}) must hold. Thus, any algorithm $\Sigma(\alpha,\beta,m,L)$ which satisfies the conditions of the theorem will have a convergence rate which satisfies (\ref{maxrateineq}).
\hfill $\blacksquare$

Also, we are now in the position to prove Theorem \ref{maintheo}.

\textit{Proof of Theorem \ref{maintheo}}: 
    In order to prove this theorem, we construct an algorithm $\Sigma(\alpha^*,\beta^*,m,L)$ satisfying the conditions of the theorem and such that (\ref{maintheorate}) is satisfied. It follows from Lemma \ref{convergespectral} that we want to construct an algorithm $\Sigma(\alpha^*,\beta^*,m,L)$ such that the characteristic polynomial of $\tilde{A_0}(\lambda)$ has all its roots inside the circle of radius $\rho = \sqrt[n]{\rho_{\text{\scriptsize HB}}} \in (0,1)$ for all $\lambda \in [m,L]$. We now define
    \begin{equation*}
        P(z) = \frac{1}{(z-1)^n}
    \end{equation*}
    and 
    \begin{equation*}
        K(z) = \frac{(z-1)^{n-1} \tilde{N}(z)}{\tilde{D}(z)}
    \end{equation*}
    where $\tilde{N}(z)$ and $\tilde{D}(z)$ are defined in (\ref{tilD}).
    It follows from (\ref{tilD}) that
    \begin{equation}
        \tilde{P}(z)\tilde{K}(z) = P(z)K(z).
    \end{equation}
    Therefore, $\sup_{m \leq \lambda \leq L} \rho(\tilde{A_0}(\lambda))$ is the radius of the smallest disk that contains the poles of the SISO feedback control system in Figure \ref{fig:BD}. 

    We now construct the interpolation function $\hat{S}(z)$ defined in (\ref{hatS=theta}) using the finite Blaschke product approach, as described in \cite[p.~53]{BoundAnalyticFunc}. The Blaschke product is a classical method for constructing bounded analytic functions on the unit disk that satisfy prescribed interpolation conditions. The general form of a finite Blaschke product is given by: 
        \begin{equation} 
            \mathbf{B}(z) = \psi \prod_{i=1}^{n} \frac{z - z_i}{1 - \overline{z_i} z}, \quad i = 1,2,...,n 
        \end{equation}
    where $|\psi|=1$, $z_i$ are the interpolation points inside the unit disk $\mathscr{D}$ and $\overline{z_i}$ denotes the complex conjugate of $z_i$. 
    Considering the interpolation conditions (\ref{modifyinterpoint1}), the corresponding Blaschke product is given by:
        \begin{equation} 
           \mathbf{B}(z) =  \psi \prod_i^n \left( \frac{z - \tfrac{\rho}{1 + \varepsilon_i}}{1 - \tfrac{\rho}{1 + \varepsilon_i} z} \right). 
        \end{equation}
    In addition, considering (\ref{modifyinterpoint2}), evaluating the Blaschke product at $z=0$, we obtain
         \begin{align} 
            \mathbf{B}(0) &= \psi \prod_i^n \left( \frac{ -\tfrac{\rho}{1 + \varepsilon_i} }{1} \right) = \psi \frac{ \rho^n }{(-1)^n\prod_i^n( 1 + \varepsilon_i) } = \theta(1). 
         \end{align}
    By letting $\varepsilon_i \rightarrow 0$, we obtain
        \begin{equation}
            \psi = (-1)^n\frac{\theta(1)}{\rho^n},
        \end{equation}
    where $|\psi| = \frac{\theta(1)}{\rho^n} =\frac{\theta(1)}{\rho_{\text{\scriptsize HB}}}= 1$. Then
        \begin{equation}
            \mathbf{B}(z) = (-1)^n\frac{\theta(1)}{\rho^n} \left(\frac{z-\rho}{1-\rho z} \right)^n =\frac{\theta(1)}{\rho^n} \left(\frac{\rho-z}{1-\rho z} \right)^n.
        \end{equation}
    Therefore, the required analytic function $\hat{S}(z)$ is 
        \begin{equation}
            \label{solutionsz}
            \hat{S}(z) =\frac{\theta(1)}{\rho^n} \left(\frac{\rho-z}{1-\rho z} \right)^n =\left(\frac{\rho-z}{1-\rho z} \right)^n ,
        \end{equation}
    since $\rho^n = \rho_{\text{\scriptsize HB}}=\theta(1)$.
    Thus, the function $\hat{S}(z)$ defined in (\ref{solutionsz}) satisfies the interpolation conditions (\ref{hatsrhointerpoint}) and (\ref{hats0interpoint}), $|\hat{S}(z)| < 1$ for all $z \in \mathscr{D}$, and inequality (\ref{rhorhohb}) is satisfied with equality.
    
    Since $\hat{S}(z) = \theta(S(\rho z^{-1}))$, we can now construct
    \begin{equation}
        \label{szinv}
        S(z) = \theta ^{-1} \left(\hat{S}(\rho z^{-1}) \right).
    \end{equation}
    In order to find $S(z)$, we firstly calculate $\hat{S}(\rho z^{-1})$:
    \begin{equation}
        \label{hatspz}
        \hat{S}(\rho z^{-1}) = \theta (1) \frac{ (z-1)^n}{\left(z-\rho^2 \right)^n}.
    \end{equation}
    Also, we construct the inverse of the function $\theta(\cdot)$ as
    \begin{equation}
        \label{thetazinv}
        \theta^{-1}(z) = \frac{8 L m z}{(L-m) \left(L (z-1)^2+m (z+1)^2\right)}.
    \end{equation}
    Substituting (\ref{hatspz}) and (\ref{thetazinv}) into (\ref{szinv}) yields
    \begin{align}
        S(z)
        &=\frac{1}{\left(L \left(\theta(1)  \left(\frac{z-1}{z-\rho ^2}\right)^n-1\right)^2+m \left(\theta(1)  \left(\frac{z-1}{z-\rho ^2}\right)^n+1\right)^2\right)} \nonumber \\
        & \quad \times\frac{8 \theta(1)  L m \left(\frac{z-1}{z-\rho ^2}\right)^n}{(L-m)}.\label{interposz}
    \end{align}
    Substituting (\ref{interposz}) and $P(z) = \frac{1}{(z-1)^n}$ into (\ref{sensifunc}) and rearranging yields
    \begin{align}
     \label{kz2}
        K(z) = \frac{(z-1)^n\left(\frac{z-1}{z-\rho ^2}\right)^{-n} K_n(z) }{4 \theta(1)  L m} 
    \end{align}
    where 
    \begin{align*}
        K_n(z)&= L \left(\theta(1) \left(\frac{z-1}{z-\rho ^2}\right)^n-1\right)^2 \\&\quad \quad-m \left(\theta(1)  \left(\frac{z-1}{z-\rho ^2}\right)^n+1\right)^2 .
    \end{align*}
    Substituting (\ref{kz2}) into $\tilde{K}(z) = \frac{1}{z-1}K(z)$ yields
    \begin{align}   
        \label{tilkz}
        \tilde{K}(z) &= \frac{1}{(z-1)^{n-1}}\frac{(z-1)^n \left(\frac{z-1}{z-\rho ^2}\right)^{-n} K_n(z)}{4 \theta(1)  L m} \nonumber\\
        &=\frac{(z-1) \left(\frac{z-\rho^2}{z-1}\right)^{n} K_n(z)}{4 \theta(1)  L m}
    \end{align}
    Substituting $\theta(1) =\rho_{\text{\scriptsize HB}} = \frac{\sqrt{L}-\sqrt{m}}{\sqrt{L}+\sqrt{m}}$ into (\ref{tilkz}) yields
    \begin{align}
    \label{tilkz1}
        \tilde{K}(z) &=\frac{(z-1) \left(\frac{z-\rho^2}{z-1}\right)^{n} \tilde{K}_n(z)}{4 L m} \nonumber \\
        &=\frac{(z-\rho^2)^n \tilde{K}_n(z)}{4Lm(z-1)^{n-1}},
    \end{align}
    where
    \begin{align}
    \label{tilknz}
        &\tilde{K}_n(z) \nonumber\\
        &=-2 (L+m) \left(\frac{z-1}{z-\rho^2}\right)^n+\left(\sqrt{L}-\sqrt{m}\right)^2 \left(\frac{z-1}{z-\rho ^2}\right)^{2 n} \nonumber \\& \quad+\left(\sqrt{L}+\sqrt{m}\right)^2 \nonumber\\
        &=\frac{1}{(z-\rho^2)^{2n}}\nonumber \\
        & \quad \times \Big( \left(\sqrt{L}+\sqrt{m}\right)^2(z-\rho^2)^{2n}+\left(\sqrt{L}-\sqrt{m}\right)^2(z-1)^{2n} \nonumber \\&\quad-2 (L+m)(z-\rho^2)^n(z-1)^n \Big).
    \end{align}
    Substituting (\ref{tilknz}) into (\ref{tilkz1}) and simplifying, yields
    \begin{align}
    \label{tilKz}
        &\tilde{K}(z)\nonumber \\&=\frac{\tilde{K}_1 (z-\rho^2)^{2n}+\tilde{K}_2 (z-1)^{2n}-\tilde{K}_3 (z-\rho^2)^n (z-1)^n}{(z-\rho^2)^n (z-1)^{n-1}},        
    \end{align}
    where $\tilde{K}_1,\tilde{K}_2,\tilde{K}_3$ are defined in (\ref{correstilK}). Note that the transfer function $\tilde{K}(z)$ is proper, since $\tilde{K}_3 = \tilde{K}_1+\tilde{K}_2$, which cancels out the leading term $z^{2n}$ in the numerator. 
    In order to find the coefficients $\alpha^*_j$ and $\beta^*_j$, we employ the binomial theorem \cite{concretemath}, which yields
    \begin{equation}
        (z-\rho^2)^{2n}=\sum_{j=0}^{2n}\binom{2n}{j}\,z^{\,2n-j}(-\rho^2)^j,
    \end{equation}
    \begin{equation}
        (z-1)^{2n}=\sum_{k=0}^{2n}\binom{2n}{k}\,z^{\,2n-k}(-1)^k,
    \end{equation}
    \begin{align}
        &(z-\rho^2)^n (z-1)^n \nonumber \\&=\sum_{i=0}^{n}\sum_{j=0}^{n}\binom{n}{i}\binom{n}{j}\,z^{\,2n-(i+j)}(-\rho^2)^i(-1)^j.
    \end{align}
    Therefore the numerator of $\tilde{K}(z)$ can be rewritten as
    \begin{align}
    \label{Numtilk}
        \tilde{N}(z) &= \tilde{K}_1\sum_{j=0}^{2n}\binom{2n}{j}z^{2n-j}(-\rho^2)^j \nonumber\\
        &\quad +\tilde{K}_2\sum_{k=0}^{2n}\binom{2n}{k}z^{2n-k}(-1)^k \nonumber \\
        &\quad -\tilde{K}_3\sum_{i=0}^{n}\sum_{j=0}^{n}\binom{n}{i}\binom{n}{j}z^{2n-(i+j)}(-\rho^2)^i(-1)^j \nonumber\\
        &= \sum_{r=0}^{2n}\Biggl( \tilde{K}_1\binom{2n}{r}(-\rho^2)^r\nonumber + \tilde{K}_2\binom{2n}{r}(-1)^r \nonumber\\ &\quad- \tilde{K}_3\sum_{i=0}^{r}\binom{n}{i}\binom{n}{r-i}(-\rho^2)^i(-1)^{r-i}\Biggr)z^{2n-r}.
    \end{align}
    Similarly,
    \begin{align}
        (z-\rho^2)^n&=\sum_{i=0}^{n}\binom{n}{i}\,z^{\,n-i}(-\rho^2)^i,\\
        (z-1)^{n-1}&=\sum_{k=0}^{n-1}\binom{n-1}{k}\,z^{\,n-1-k}(-1)^k,
    \end{align}
    and the denominator of $\tilde{K}(z)$ can be rewritten as 
    \begin{equation}
    \label{Dentildk}
        \tilde{D}(z) = \sum_{i=0}^{n}\sum_{k=0}^{n-1}\binom{n}{i}\binom{n-1}{k}\,z^{\,2n-1-(i+k)}(-\rho^2)^i(-1)^k.
    \end{equation}
    Note that 
    \begin{equation}
        \rho = \theta(1)^{\frac{1}{n}} = \left(\frac{\sqrt{L}-\sqrt{m}}{\sqrt{L}+\sqrt{m}} \right)^{\frac{1}{n}}.
    \end{equation}
    Therefore, we obtain the final expression for $\tilde{K}(z)$
    \begin{equation}
        \label{tilkzfinal}
        \tilde{K}(z) = \frac{\tilde{N}(z)}{\tilde{D}(z)},
    \end{equation}
    where the coefficients $\alpha_j$ and $\beta_j$ in the polynomial expansions of $\tilde{N}(z)$ and $\tilde{D}(z)$ precisely match the corresponding coefficients obtained from the binomial expansions of $\tilde{N}(z)$ and $\tilde{D}(z)$ in (\ref{Numtilk}) and (\ref{Dentildk}), respectively.
    Thus, we have constructed an algorithm $\Sigma(\alpha^*,\beta^*,m,L)$ satisfying the conditions of the theorem with convergence rate given by (\ref{maintheorate}). \hfill $\blacksquare$

\begin{remark}
    Note that it may be useful to define an alternative form for the algorithm (\ref{uprule}) by using an alternative state space realization 
    \begin{align*}
        &\tilde{C}_0(zI-A_0)^{-1}\tilde{B}_0\\ &= \tilde{P}(z)\tilde{K}(z)\\
        &=\frac{\tilde{K}_1 (z-\rho^2)^{2n}+\tilde{K}_2 (z-1)^{2n}-\tilde{K}_3 (z-\rho^2)^n (z-1)^n}{(z-\rho^2)^n (z-1)^{n}}
    \end{align*}
    rather than using the formula (\ref{eq:alpha}), (\ref{eq:beta}). This is an area for future research.
\end{remark}


\section{Illustrative Example}
\label{IlluExam}
This section gives an illustrative example of the use of the algorithm $\Sigma(\alpha^*,\beta^*,m,L)$ of the form (\ref{uprule}) defined in Theorem \ref{maintheo}. Consider a quadratic cost function $f(x,t) \in \mathscr{F}^q_{m,L}$,
\begin{equation}
    f(x,t) = \frac{1}{2}(x(t)-x^*(t))^2
\end{equation}
where $x^*(t)\in \mathbb{R}$ is a real source location at time $t$ defined by
\begin{equation}
    x^*(t) = \tanh{t} = \frac{e^t-e^{-t}}{e^t+e^{-t}}. \label{tanhxt}
\end{equation}
Equation (\ref{tanhxt}) has a Taylor series expansion \cite[Chapter 4]{HandbookMF}
\begin{align}
    \tanh{t}
    &= t-\frac{1}{3}t^3+\frac{2}{15}t^5-\frac{17}{315}t^7+... .\label{TSpoly}
\end{align}
To simplify the problem, we only consider the first two terms of the Taylor series (\ref{TSpoly}). It follows from Lemma \ref{maintheonintegrators} that to achieve exact tracking of a polynomially varying point of order $3$, the iterative method must incorporate at least four integrators. In order to find the corresponding algorithm parameters, it follows from (\ref{tilD}) and (\ref{tilKz}) that
\begin{align}
     \sum_{j=0}^{k}\alpha_j z^{k-j}&=\tilde{K}_1 (z-\rho^2)^{2n}+\tilde{K}_2 (z-1)^{2n} \nonumber\\&\quad-\tilde{K}_3 (z-\rho^2)^n (z-1)^n,\\
     z^k - \sum_{j=0}^{k-1} \beta_j z^{k-j-1}&=(z-\rho^2)^n (z-1)^{n-1}
\end{align}
where $\tilde{K}_1, \tilde{K}_2, \tilde{K}_3$ are defined in (\ref{correstilK}) and  $n=4$. Solving these equations, we find the corresponding optimal parameters as follows:
\begin{align*}
    \beta^*_0 &= 4\rho^2 + 3, \\
    \beta^*_1 &= -6\rho^4 - 12\rho^2 - 3, \\
    \beta^*_2 &= 4\rho^6 + 18\rho^4 + 12\rho^2 + 1, \\
    \beta^*_3 &= -\rho^8 - 12\rho^6 - 18\rho^4 - 4\rho^2, \\
    \beta^*_4 &= 3\rho^8 + 12\rho^6 + 6\rho^4, \\
    \beta^*_5 &= -3\rho^8 - 4\rho^6, \\
    \beta^*_6 &= \rho^8,
\end{align*}
and
\begin{align*}
    \alpha^*_0 &= \frac{4}{\sqrt{L}\sqrt{m}} - \frac{4\rho^2}{\sqrt{L}\sqrt{m}}, \\
    \alpha^*_1 &= \frac{14\rho^4 - 14}{\sqrt{L}\sqrt{m}} 
    + \frac{4\rho^4 - 8\rho^2 + 4}{L} 
    + \frac{4\rho^4 - 8\rho^2 + 4}{m}, \\
    \alpha^*_2 &= \frac{-28\rho^6 + 28}{\sqrt{L}\sqrt{m}} 
    + \frac{-12\rho^6 + 12\rho^4 + 12\rho^2 - 12}{L} 
    \\ &\quad + \frac{-12\rho^6 + 12\rho^4 + 12\rho^2 - 12}{m}, \\
    \alpha^*_3 &= \frac{35\rho^8 - 35}{\sqrt{L}\sqrt{m}} 
    + \frac{17\rho^8 - 8\rho^6 - 18\rho^4 - 8\rho^2 + 17}{L} 
    \\ &\quad+ \frac{17\rho^8 - 8\rho^6 - 18\rho^4 - 8\rho^2 + 17}{m}, \\
    \alpha^*_4 &= \frac{-28\rho^{10} + 28}{\sqrt{L}\sqrt{m}} 
    \\&\quad+ \frac{-14\rho^{10} + 2\rho^8 + 12\rho^6 + 12\rho^4 + 2\rho^2 - 14}{L} \\
    &\quad + \frac{-14\rho^{10} + 2\rho^8 + 12\rho^6 + 12\rho^4 + 2\rho^2 - 14}{m}, \\
    \alpha^*_5 &= \frac{14\rho^{12} - 14}{\sqrt{L}\sqrt{m}} 
    + \frac{7\rho^{12} - 3\rho^8 - 8\rho^6 - 3\rho^4 + 7}{L} 
    \\&\quad+ \frac{7\rho^{12} - 3\rho^8 - 8\rho^6 - 3\rho^4 + 7}{m}, \\
    \alpha^*_6 &= \frac{-4\rho^{14} + 4}{\sqrt{L}\sqrt{m}} 
    + \frac{-2\rho^{14} + 2\rho^8 + 2\rho^6 - 2}{L} 
    \\&\quad+ \frac{-2\rho^{14} + 2\rho^8 + 2\rho^6 - 2}{m}, \\
    \alpha^*_7 &= \frac{\rho^{16} - 1}{2\sqrt{L}\sqrt{m}} 
    + \frac{\rho^{16} - 2\rho^8 + 1}{4L} 
    + \frac{\rho^{16} - 2\rho^8 + 1}{4m},
\end{align*}
where 
\begin{equation*}
    \rho = \left(\frac{\sqrt{L} - \sqrt{m}}{\sqrt{L} + \sqrt{m}}\right)^{\frac{1}{4}}.
\end{equation*}

\begin{figure}[htb!]
    \centering
    \includegraphics[width=\linewidth]{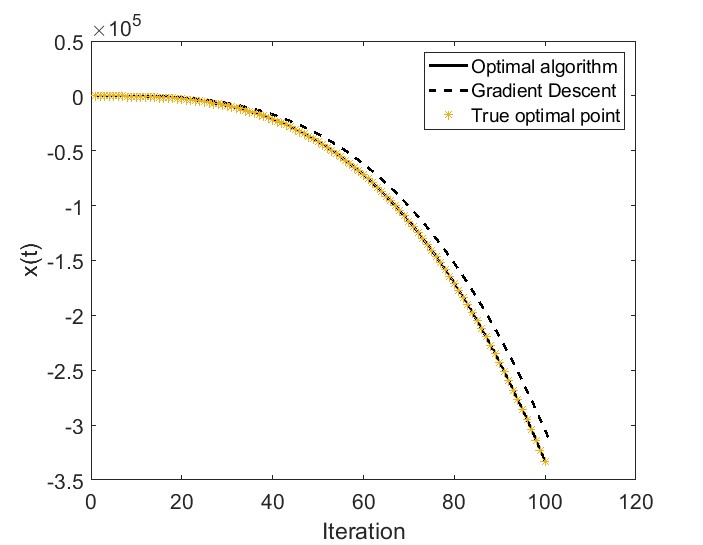}
    \caption{Comparison between the true optimal point and its estimated value using the optimal algorithm from Theorem \ref{maintheo} and the gradient descent method \cite{Boyd_Vandenberghe_2004}, where $L=5$ and $m=1$.}
    \label{fig:track3}
\end{figure}
Figure \ref{fig:track3} compares the true optimal point to its estimated value obtained by using Theorem \ref{maintheo} and by using the standard gradient descent method with step size $\frac{2}{L+m}$, where $L=1$ and $m=5$.  
The yellow star represents the true optimal point, moving according to the polynomial trajectory (\ref{TSpoly}) of order $3$. The red line shows the optimal point estimated by the optimal method of Theorem \ref{maintheo}. It can be observed that the optimal method of Theorem \ref{maintheo} successfully tracks the optimal point with zero steady-state error. In addition, the black line in Figure \ref{fig:track3} shows the performance of the gradient descent method, which fails to track the time-varying optimal point with zero steady-state error.

\begin{figure}[htb!]
    \centering
    \includegraphics[width=\linewidth]{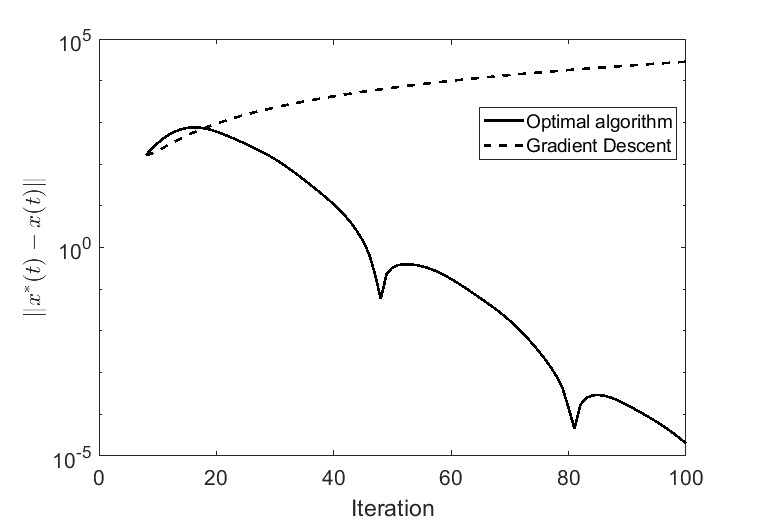}
    \caption{
    Difference between the true optimal point and its estimate using the optimal method of Theorem \ref{maintheo} and the gradient descent method \cite{Boyd_Vandenberghe_2004}.}
    \label{fig:error3}
\end{figure}
To further verify the performance of our approach, Figure \ref{fig:error3} illustrates the difference between the true optimal point, its estimate using the optimal algorithm obtained from Theorem \ref{maintheo} and the gradient descent method. Clearly, the optimal algorithm obtained from Theorem \ref{maintheo} achieves zero steady-state error, whereas the gradient descent method does not.

\section{CONCLUSIONS}\label{sec5}
In this paper, we establish a fundamental convergence rate bound for all linear gradient based optimization algorithms which optimize a time-varying cost function with zero steady-state error. We consider online optimization algorithms for solving optimization problems with time-varying quadratic cost functions. Our main results incorporate results on the optimal gain margin of linear uncertain systems. We also use a version of the internal model principle to require the inclusion of discrete time integrators in the corresponding algorithm. This approach enables the algorithm to achieve exact tracking of time-varying optimal points. 









\bibliographystyle{IEEEtran}
\bibliography{IEEEabrv,main}

\begin{thebibliography}{10}
\providecommand{\url}[1]{#1}
\csname url@rmstyle\endcsname
\providecommand{\newblock}{\relax}
\providecommand{\bibinfo}[2]{#2}
\providecommand\BIBentrySTDinterwordspacing{\spaceskip=0pt\relax}
\providecommand\BIBentryALTinterwordstretchfactor{4}
\providecommand\BIBentryALTinterwordspacing{\spaceskip=\fontdimen2\font plus
\BIBentryALTinterwordstretchfactor\fontdimen3\font minus \fontdimen4\font\relax}
\providecommand\BIBforeignlanguage[2]{{%
\expandafter\ifx\csname l@#1\endcsname\relax
\typeout{** WARNING: IEEEtran.bst: No hyphenation pattern has been}%
\typeout{** loaded for the language `#1'. Using the pattern for}%
\typeout{** the default language instead.}%
\else
\language=\csname l@#1\endcsname
\fi
#2}}

\bibitem{Introonline}
E.~Hazan, \emph{Introduction to Online Convex Optimization}, 2nd~ed.\hskip 1em plus 0.5em minus 0.4em\relax USA: The MIT Press, 2022.

\bibitem{Zampieri2023}
N.~Bastianello, R.~Carli, and S.~Zampieri, ``Internal model-based online optimization,'' \emph{IEEE Transactions on Automatic Control}, vol.~69, no.~1, pp. 689--696, 2024.

\bibitem{madden2021bounds}
L.~Madden, S.~Becker, and E.~Dall’Anese, ``Bounds for the tracking error of first-order online optimization methods,'' \emph{Journal of Optimization Theory and Applications}, vol. 189, no.~2, pp. 437--457, 2021.

\bibitem{doi:10.1137/15M1009597}
L.~Lessard, B.~Recht, and A.~Packard, ``Analysis and design of optimization algorithms via integral quadratic constraints,'' \emph{SIAM Journal on Optimization}, vol.~26, no.~1, pp. 57--95, 2016.

\bibitem{IQCScherer}
C.~S. Simon~Michalowsky and C.~Ebenbauer, ``Robust and structure exploiting optimisation algorithms: an integral quadratic constraint approach,'' \emph{International Journal of Control}, vol.~94, no.~11, pp. 2956--2979, 2021.

\bibitem{Lessard2023L}
B.~Van~Scoy and L.~Lessard, ``A tutorial on a {L}yapunov-based approach to the analysis of iterative optimization algorithms,'' in \emph{2023 62nd IEEE Conference on Decision and Control (CDC)}, 2023, pp. 3003--3008.

\bibitem{lessard2022}
L.~Lessard, ``The analysis of optimization algorithms, a dissipativity approach,'' \emph{IEEE Control Systems Magazine}, vol.~42, 05 2022.

\bibitem{UGRINOVSKII2023111129}
V.~Ugrinovskii, I.~R. Petersen, and I.~Shames, ``A robust control approach to asymptotic optimality of the heavy ball method for optimization of quadratic functions,'' \emph{Automatica}, vol. 155, pp. 111--129, 2023.

\bibitem{Jovanovic2024Tann}
W.~Wu, J.~Chen, M.~R. Jovanovi{\'c}, and T.~T. Georgiou, ``Tannenbaum's gain-margin optimization meets {P}olyak's heavy-ball algorithm,'' \emph{arXiv preprint arXiv:2409.19882}, 2024.

\bibitem{frequencyTannen2024}
S.~Zhang, W.~Wu, Z.~Li, J.~Chen, and T.~T. Georgiou, ``Frequency-domain analysis of distributed optimization: Fundamental convergence rate and optimal algorithm synthesis,'' \emph{IEEE Transactions on Automatic Control}, pp. 1--16, 2024.

\bibitem{Tannenbaum1985}
P.~Khargonekar and A.~Tannenbaum, ``Non-euclidian metrics and the robust stabilization of systems with parameter uncertainty,'' \emph{IEEE Transactions on Automatic Control}, vol.~30, no.~10, pp. 1005--1013, 1985.

\bibitem{AlexACCTracking}
A.~Wu, I.~Petersen, V.~Ugrinovskii, and I.~Shames, ``An online optimization algorithm for tracking a linearly varying optimal point with zero steady-state error,'' in \emph{2024 American Control Conference (ACC)}, 2025.

\bibitem{AlexECC2025}
A.~Wu, I.~Petersen, and I.~Shames, ``Online optimization with integral action: An optimal algorithm for a time varying quadratic cost function,'' in \emph{2025 European Control Conference (ECC)}, 2025.

\bibitem{concretemath}
R.~L. Graham, D.~E. Knuth, and O.~Patashnik, \emph{Concrete Mathematics: A Foundation for Computer Science}, 2nd~ed.\hskip 1em plus 0.5em minus 0.4em\relax USA: Addison-Wesley Longman Publishing Co., Inc., 1994.

\bibitem{pick1915}
G.~Pick, ``{\"U}ber die beschr{\"a}nkungen analytischer funktionen, welche durch vorgegebene funktionswerte bewirkt werden,'' \emph{Mathematische Annalen}, vol.~77, no.~1, pp. 7--23, 1915.

\bibitem{Nevanlinna1919}
R.~Nevanlinna, ``Über beschränkte funktionen die in gegebenen punkten vorgeschriebene werte annehmen,'' \emph{Ann. Acad. Sci. Fenn. Ser. A 1 Mat. Dissertationes}, vol.~13, 1919.

\bibitem{POLYAK19641}
B.~T. Polyak, ``Some methods of speeding up the convergence of iteration methods,'' \emph{USSR Computational Mathematics and Mathematical Physics}, vol.~4, no.~5, pp. 1--17, 1964.

\bibitem{LecOnCOptNes2018}
Y.~Nesterov, \emph{Lectures on Convex Optimization}, 2nd~ed.\hskip 1em plus 0.5em minus 0.4em\relax Springer Publishing Company, Incorporated, 2018.

\bibitem{IterativeSolution}
J.~M. Ortega and W.~C. Rheinboldt, \emph{Iterative Solution of Nonlinear Equations in Several Variables}.\hskip 1em plus 0.5em minus 0.4em\relax SIAM, 2000.

\bibitem{fadali2009digital}
M.~Fadali and A.~Visioli, \emph{Digital Control Engineering: Analysis and Design}.\hskip 1em plus 0.5em minus 0.4em\relax Academic Press, 2009.

\bibitem{BoundAnalyticFunc}
J.~B. Garnett, \emph{Bounded analytic functions}.\hskip 1em plus 0.5em minus 0.4em\relax Springer, 2011.

\bibitem{Schur1917}
J.~Schur, ``Über potenzreihen, die im innern des einheitskreises beschränkt sind.'' \emph{Journal für die reine und angewandte Mathematik}, vol. 147, pp. 205--232, 1917.

\bibitem{gendet}
T.~Amdeberhan and D.~Zeilberger, ``{"Trivializing"} generalizations of some {Izergin-Korepin}-type determinants,'' \emph{Discrete Mathematics \& Theoretical Computer Science}, vol. Vol. 9 no. 1, Jan 2007.

\bibitem{HandbookMF}
M.~Abramowitz and I.~A. Stegun, \emph{Handbook of Mathematical Functions, With Formulas, Graphs, and Mathematical Tables,}.\hskip 1em plus 0.5em minus 0.4em\relax USA: Dover Publications, Inc., 1974.

\bibitem{Boyd_Vandenberghe_2004}
S.~Boyd and L.~Vandenberghe, \emph{Convex Optimization}.\hskip 1em plus 0.5em minus 0.4em\relax Cambridge University Press, 2004.

\end{thebibliography}

\end{document}